\documentclass[11pt]{article}
\usepackage{amssymb, amsmath,  amsfonts,dsfont, mathrsfs,euscript,colonequals}
\usepackage{textcomp}
\usepackage{latexsym}
\usepackage{indentfirst}
\usepackage[normalem]{ulem}
\usepackage{graphicx}
\usepackage[colorlinks=true, allcolors=blue]{hyperref}

\newtheorem{Theorem}{Theorem}
\newtheorem{Lemma}{Lemma}

\newtheorem{Corollary}{Corollary}

\newtheorem{Definition}{Definition}

\newcommand{\N}{\mathbb N}
\newcommand{\R}{\mathbb R}
\newcommand{\Z}{\mathbb Z}

\newcommand{\CC}{\mathbb C}
\newcommand{\PP}{\text{Pr}}
\newcommand{\E}{\mathbb E}
\newcommand{\disteq}{\stackrel{d}{=} }
\newcommand{\disto}{\stackrel{d}{\longrightarrow} }
\newcommand{\defeq}{\colonequals}

\begin{document}
	\author{I. A. Alexeev \footnote{Institute for Information Transmission Problems of Russian Academy of Sciences, Bolshoy Karetny per. 19, build.1, Moscow and St. Petersburg Department 
of Steklov Mathematical Institute of Russian Academy of Sciences, 27 Fontanka, St. Petersburg, Russia, email: \texttt{vanyalexeev@list.ru}}, I.A. Melnikov\footnote{Institute for Information Transmission Problems of Russian Academy of Sciences, Bolshoy Karetny per. 19, build.1, Moscow, email: \texttt{mi@melnikov-ignat.ru}}, A.Y.Uglovski\footnote{Institute for Information Transmission Problems of Russian Academy of Sciences, Bolshoy Karetny per. 19, build.1, Moscow, email: \texttt{uglovskij.au@gmail.com}}}

	\title{Limit theorems for pseudo sum of discrete distributions.\thanks{The work was supported by the Russian Science Foundation (grant No. 22-21-00016).}}
    \date{}
	\maketitle

 \begin{center}
     \textbf{Annotation}
 \end{center}
In this article we introduce associative Look-Up Tables. With their help, pseudo sums are correctly determined. The set of limit distributions in a pseudo-summation scheme of i.i.d. random variables is described. Also, two special cases that are similar to the classical sum and maximum operations are considered. In both situations, the set of stable distributions and the set of infinity divisible distributions are fully described. In addition, necessary and sufficient conditions for random variable to belong to the domain of attraction of stable random variable are introduced.

	\section{Introduction}
The present paper is devoted to the study of discrete stable distributions in the case when classical summation is changed to some associative operation. 

In the early 1930s P. L\'evy described the set of limit distributions in a summation scheme of independent identically distributed (i.i.d.) random variables with some positive normalization and real-valued centering (see \cite{Levy}). He proved that if for some $B_n >0$ and $a_n \in \R$, $n \in \N$ we have
\begin{equation*}
    \frac1{B_n} \sum_{k=1}^n X_k - a_n \disto \xi,  
\end{equation*}
then $\xi$ has to be stable. Here and after, $\disto$ denotes the convergence in distribution. 

Recall that, a random variable $\xi$ is called stable if for every $b_1$, $b_2 > 0$ there exist $b >0$ and $a \in \R$ such that 
\begin{equation*}
    b_1 \xi_1 + b_2 \xi_2 \disteq b\xi + a,
\end{equation*}
where $\xi_1$, $\xi_2$ are independent copies of $\xi$ and $\disteq $ denotes the equality in distribution. The class of classical one-dimensional stable distributions is well studied.

It has been noted that the sum operation does not play such a significant role in L\'evi's results. To a greater extent, an associative, commutative and reversible operation is required. For example, similar results exist for maximum stable laws (see \cite{Resnik}). More precisely, if for some $B_n >0$, $n \in \N$ we have
\begin{equation*}
    \frac1{B_n} \max_{k=1,\ldots,n} X_k \disto \xi,  
\end{equation*}
then $\xi$ has to be max-stable. 

Both stable and max-stable distributions are well studied and fully characterized (see \cite{Gnedenko}, \cite{Ibragimov}, \cite{Levy}, \cite{Meer}, \cite{Resnik}, \cite{Sato}). 

The motivation for this article was the results on stable laws in all their understandings and the theory of Low-Density Parity Check codes. Namely, in standard Information Bottleneck decoder so called Look-Up Tables are used. Firstly, quantize data is used when transmitting a message, that is, only a finite number of values $\chi = \{x_1, \ldots, x_N\}$ are transmitted and some special sum operation $\oplus$ is used when updating messages. This article is devoted to such sum operations and limit theorems with them. 

Similar limit theorems have already been studied before. In the case when the operation is associative, then the set $(\chi, \oplus)$ with the operation is a semigroup. Limit theorems with a semigroup operation were studied in the \cite{Lof}. The case of a non-associative operation was studied in the articles \cite{Markovski}, \cite{Yashunsky} from the point of view of a Markov chain with a finite number of values.

In section 2 associative Look-Up Tables are introduced and the most general properties are proved. In particular, the set of limit distributions in a pseudo-summation scheme of i.i.d. random variables is described. Section 3 is devoted to two special cases that are similar to the classical sum and maximum operations. In both situations, the set of stable distributions and the set of infinity divisible distributions are fully described. Also, necessary and sufficient conditions for random variable to belong to the domain of attraction of stable random variable are introduced. 

\section{General results}
Let $\chi = \{ x_1, \ldots, x_N\}$, $N \in \N$ be a finite set of distinct real values, i.e. $x_j \in \R$, $j = 1,\ldots, N$, $x_j \ne x_k$ for $j \ne k$. Every function $f\colon \chi\times \chi \to \chi$ can be presented as a matrix  $A = (A(j,k))_{j,k=1}^N \in \chi^{N \times N}$, where $f(x_i, x_j) = A(i,j)$ for all $i,j = 1,\ldots, N$. 

Let assume that for all $i,j,k = 1,\ldots, N$ the following condition holds:
\begin{equation}\label{eq:associativity}
    A(i, A(j,k)) = A(A(i,j), k).
\end{equation}

It follows that function $f$ (or matrix $A$) can correctly denote a pseudo-summation operation:
\begin{equation*}
    x_j \oplus x_k \defeq A(j,k) = f(x_i, x_j).
\end{equation*} 

Condition \ref{eq:associativity} yields the associativity of a pseudo-summation, that is 
\begin{equation*}
    x \oplus \bigr( y \oplus z \bigr) = \bigr( x \oplus y\bigr) \oplus z, \quad x,y,z \in \chi.
\end{equation*}

Let us note that such an operation is commonly used in the coding theory, namely, in the standard Information Bottleneck decoder for Low-Density Parity Check codes (for further information see, for instance, \cite{He2019}). The matrix $A$ is called \textit{Look-Up Table}. Further, this is the name that will be used. Also, if for $A$ condition \ref{eq:associativity} holds, then $A$ will be called \textit{associative Look-Up Table}. 

For each associative Look-Up Table, similarly to the classical real case, stable random variables can be determined for the pseudo-summation.
\begin{Definition}
A random variable $\xi$ taking values in $\chi$ is called stable if $\xi_1 \oplus \xi_2 \disteq \xi$, where $\xi_1$, $\xi_2$ are independent copies of $\xi$. 
\end{Definition}

One can easily show that $\xi$ is stable if for every $m \in \N$ the sum $\xi_1 \oplus \ldots \oplus \xi_m \disteq \xi$, where $\xi_1, \ldots, \xi_m$ are independent copies of $\xi$. 

The set of stable distributions directly depends on the Look-Up Table $A$. Thus, we will write that random variable $\xi$ is $A$-stable in order to highlight with respect to which operation the random variable is stable.

Similar to the classic sum, stable distributions pays the crucial role in the limit theorems for pseudo-summation.
\begin{Theorem}
    Let $\{X_k\}$, $k \in \N$ be a sequence of i.i.d. random variables taking values in $\chi$ and let $A$ be an associative Look-Up Table. If for some random variable $\xi$ we have
    \begin{equation*}
        X_1 \oplus \ldots \oplus X_m \disto \xi, \quad m \to \infty, 
    \end{equation*}
    then $\xi$ is $A$-stable. 
\end{Theorem}
\textbf{Proof.}\quad The proof is almost no different from the real case. Since $A$ is an associative Look-Up Table, then 
\begin{equation}\label{eq:two_sums}
    X_1 \oplus \ldots \oplus X_m = \bigr( X_1 \oplus \ldots \oplus X_l\bigr) \oplus \bigr( X_{l+1} \oplus \ldots \oplus X_m\bigr) \disto \xi_1 \oplus \xi_2,
\end{equation}
where $l = [m/2]$, $\xi_1$ and $\xi_2$ are independent copies of $\xi$. Since the left hand-side of \eqref{eq:two_sums} converges to $\xi$, then $\xi_1 \oplus \xi_2 \disteq \xi$, which means that $\xi$ is $A$-stable. $\Box$

Let assume that there exists $J \subset \chi$ such that for every $a,b \in J$ the equation $x \oplus a = b$ has a unique solution $x \in J$. In other words, on some subset $J$ a left pseudo-subtraction can be defined. If $\xi_1$, $\xi_2$ are independent random variables and uniformly distributed on the set $J$, then for all $b \in J$ we have 
\begin{equation*}
    \PP\bigr(\xi_1 \oplus \xi_2 = b \bigr) = \sum_{a \in J} \PP(\xi_1 \oplus a = b)\PP(\xi_2 = a) = \frac1{|J|}\sum_{a \in J}\PP(\xi_1 \oplus a = b), 
\end{equation*}
where $|J|$ is the number of elements in $J$.

Since for every $a,b \in J$ there exists unique $x \in J$ such that $x \oplus a = b$, then $\PP(\xi_1 \oplus a = b) = \frac1{|J|}$ and, hence, $\xi_1 \oplus \xi_2$ is also a uniformly distributed random variable. It follows that if on some subset $J$ a left pseudo-subtraction can be defined, then random variable that uniformly distributed on $J$ is $A$-stable. 

Let assume that $J = \{x\}$, $x \in \chi$. Then if $x \oplus x = x$, it follows that the degenerate law at point $x$ is $A$-stable.

On the other hand, if degenerate law at point $x$ is $A$-stable, then 
\begin{equation*}
    1 = \PP\bigr( \xi_1 \oplus \xi_2 = x\bigr) = \begin{cases} 0, & x \oplus x \ne x; \\ 1, & x \oplus x = x \end{cases}. 
\end{equation*}
It follows that $x \oplus x = x$ and, consequently, the degenerate law at point $x$ is $A$-stable if and only if $x \oplus x = x$.

In addition to associativity, we introduce two more conditions. Let assume that $\chi = [0:N] = [0, N) \cap \Z = \{0,\ldots,  N-1\}$ and that for an associative Look-Up Table $A$ there exist two matrices $U = (u_{i,j}) \in \CC^{N \times N}$ and $V = (v_{i,j})\in \CC^{N \times N}$ such that $\det(U) = 1$, $|v_{i,j}| \leqslant 1$ for every $i,j = 0, \ldots, N-1$, and for every $i,j$ we have
\begin{equation}\label{eq:general_ch_fun}
    u_{i \oplus j, k} = v_{i,k} \cdot u_{j,k}; \quad v_{i \oplus j, k} = v_{i,k} \cdot v_{j,k}, \,\, k = 0,\ldots,N-1, 
\end{equation}
where $\oplus$ is an operation generated by $A$. 

\begin{Lemma}\label{lemma:ch_fun}
If Look-Up Table $A$ satisfy \eqref{eq:general_ch_fun}, then for every independent random variables $X$, $Y$ taking values in $\chi$ we have
\begin{equation}\label{eq:gen_ch_fun}
    U\bigr(\mathcal{L}(X\oplus Y)\bigr) = V\bigr(\mathcal{L}(X)\bigr) \odot U\bigr(\mathcal{L}(Y)\bigr), \quad V\bigr(\mathcal{L}(X\oplus Y)\bigr) = V\bigr(\mathcal{L}(X)\bigr) \odot V\bigr(\mathcal{L}(Y)\bigr),
\end{equation}
where $\mathcal{L}(X) = \bigr(\PP(X = i)\bigr)$, $i = 0, \ldots, N-1$ is a vector in $\R^N$ corresponding to the distribution of $X$, by $\odot$ we denote the component-by-component multiplication.  
\end{Lemma}
\textbf{Proof.} \quad For simplicity, let $\mathcal{L}(X) = P = (p_0, \ldots, p_{N-1})$, $\mathcal{L}(Y) = Q = (q_0, \ldots, q_{N-1})$, and $\mathcal{L}(X \oplus Y) = R = (r_0, \ldots, r_{N-1})$. Then $P = \sum_{k=0}^{N-1}p_k e_k$, $Q = \sum_{j=0}^{N-1}q_j e_j$, and $R = \sum_{j=0}^{N-1}r_j e_j$, where $e_0, \ldots, e_{N-1}$ denote a standard basis in $\R^N$. Then 
\begin{equation*}
    V(P) \odot V(Q) = \Bigr( \sum_{k=0}^{N-1}p_k v_k \Bigr) \odot \Bigr( \sum_{j=0}^{N-1}q_j v_j \Bigr) = \sum_{k,j = 0}^{N-1} p_k q_j v_k \odot v_j = \sum_{k,j = 0}^{N-1} p_k q_j v_{k\oplus j} = V(R),
\end{equation*}
where $v_j$ is a $j$-th column of a matrix $V$. 

Similarly one can show that $U(R) = V(P) \odot U(Q)$. $\Box$

Lemma \eqref{lemma:ch_fun} shows that if Look-Up Table $A$ satisfy condition \eqref{eq:general_ch_fun}, then there exists a characteristic function. Indeed, if $P \in \R^N$ is a distribution of random variable $X$, then $U(P)$ uniquely defines $P$ and for any other distribution $Q \in \R^N$ we have
\begin{equation*}
    U(P*Q) = V(P) \odot U(Q); \quad V(P*Q) = V(P) \odot V(Q),
\end{equation*}
where $P*Q$ is a distribution of pseudo sum of independent random variables with distributions $P$ and $Q$ respectively. In other words, such a pair of matrices diagonalizes the pseudo-convolution operator. 

Let us note that not all Look-Up Tables satisfy \eqref{eq:general_ch_fun} and accordingly \eqref{eq:gen_ch_fun}. For instance, if $i \oplus j = i$ for all $i,j \in \chi$, then $v_{i,l} = v_{i \oplus i, l} = v_{i,l}v_{i,l} \in \{0,1\}$ for all $i,j,l = 0,\ldots, N-1$. Since $i \oplus j = i$, then for all independent random variables $X$, $Y$ we have $\mathcal{L}(X\oplus Y) = \mathcal{L}(Y)$ and 
\begin{equation*}
    U\bigr(\mathcal{L}(X)\bigr) = V\bigr(\mathcal{L}(X)\bigr) \odot U\bigr(\mathcal{L}(Y)\bigr), \quad V\bigr(\mathcal{L}(X)\bigr) = V\bigr(\mathcal{L}(X)\bigr) \odot V\bigr(\mathcal{L}(Y)\bigr).
\end{equation*}

Then for any random variable $X$ the vector $V\bigr(\mathcal{L}(X)\bigr) \in \{0,1\}^N$. It follows that either $v_{i,l}  = 1$ or $v_{i,l}  = 0$ for all $i,l = 0,\ldots, N-1$. The latter instantly leads to a contradiction with the invertibility of the matrix $U$. Then we have that for all independent random variables $X$ and $Y$ we have $U\bigr(\mathcal{L}(X)\bigr) = U\bigr(\mathcal{L}(Y)\bigr)$. $U$ is an invertable matrix, then $X \disteq Y$, which leads to a contradiction. 

It is trivial that $X_m \disto X$ if and only if $U(P_m) \to U(P)$, where $P_m$, $P$ are the distributions of $X_m$, $X$ respectively. Also, if $X_m \to X$, then $V(P_m) \to V(P)$.  

The following definition is a trivial generalization of the classical one.
\begin{Definition}
    A random variable $X$ taking values in $\chi$ is said to belong to the domain of attraction of a random variable $\xi$ if 
\begin{equation*}
    X_1 \oplus \ldots \oplus X_m \disto \xi. 
\end{equation*}
\end{Definition}

The following theorem shows the necessary condition for random variable $X$ to belong to the domain of attraction of another random variable $\xi$.
\begin{Theorem}\label{th:dom_attr}
    Let $A$ be associative Look-Up Table satisfying \eqref{eq:general_ch_fun}. If a random variable $X$ with distribution $P$ belongs to the domain of attraction of the random variable $\xi$ with distribution $\Xi$, then $\Xi*P = \Xi$, i.e. $\xi \oplus X \disteq \xi$. 
\end{Theorem}
\textbf{Proof.} \quad Let $X$ belongs to the domain of attraction of random variable $\xi$. Then 
\begin{equation*}
    U(P_1 * \ldots * P_m) = V(P_1 * \ldots * P_{m-1}) \odot U(P) =  \odot \bigr(V(P)\bigr)^{\odot(m-1)} \odot U(P) \to U(\Xi)\,\, \text{and}\,\, \bigr(V(P)\bigr)^{\odot m} \to V(\Xi).
\end{equation*}

It follows that $U(\Xi*P) = V(\Xi) \odot U(P) = U(\Xi)$. Since $U(P)$ defines $P$ uniquely, then $\Xi*P = P$ or $\xi \oplus X \disteq \xi$. $\Box$

\begin{Corollary}
    Let $A$ be associative Look-Up Table satisfying \eqref{eq:general_ch_fun} and $x\oplus x = x$ for some $x \in \chi$. If a random variable $X$ belongs to the domain of attraction of the degenerate law at point $\xi$, then $\PP( x \oplus X = x) = 1$.
\end{Corollary}
\section{Special cases}
In this paragraph two special cases will be considered. 
As previously, let us assume that  $\chi = [0:N]$. 

\textbf{A)} Firstly, let us consider the following Look-Up Table: 
\begin{equation}\label{eq:main_lut}
    x \oplus y = A(x,y) \defeq s^{-1} \Bigr( \bigr(s(x) + s(y) \bigr) \% N \Bigr), \quad x,y \in \chi,
\end{equation}
where $s$ is some permutation, $\%$ is a remainder of the division. 

\begin{Theorem}\label{th:operation}
    The operation $x \oplus y$ given by the formula \eqref{eq:main_lut} satisfy four following conditions:
    \begin{enumerate}
    \item[(a)] Associativity -- $A(x, A(y,z)) = A(A(x,y), z)$ for all $x,y,z \in \chi$;
    \item[(b)] Commutativity -- $A(x,y) = A(y,x)$ for all $x,y \in \chi$; 
    \item[(c)] For all $b,c \in \chi$ the equation $x \oplus b = c$ has a unique solution $x = c \ominus b$.
    \item[(d)] There exists $x \in \chi$ such that $x \oplus y = y$ for all $y \in \{0,1\ldots, N-1\}$. 
\end{enumerate}
\end{Theorem}
\textbf{Proof.} \quad First, let us show that it is sufficient to prove that  for all $x,y,z \in \chi$ we have
\begin{equation}\label{eq:group_bijective}
    x \oplus y = z \,\, \text{if and only if} \,\, e^{i\tfrac{2\pi s(x)}{N}}\cdot e^{i\tfrac{2\pi s(y)}{N}} = e^{i\tfrac{2\pi s(z)}{N}}. 
\end{equation}

If \eqref{eq:group_bijective} holds, then $\Bigr( \chi, \oplus\Bigr)$ is a group isomorphic to the $\Bigr( \{e^{i\tfrac{2\pi k}{N}}\}_{k=1}^{N-1}, \cdot\Bigr)$. It means that $\oplus$, as the $\cdot$, satisfy all four conditions $(a)-(d)$. 

Let us now prove \eqref{eq:group_bijective}. From \eqref{eq:main_lut} it follows that $x \oplus y = z$ is equivalent to $s(z) = \bigr(s(x) + s(y) \bigr) \% N$. It means that $s(z)$ is either $s(x) + s(y)$ or $s(x) + s(y) - N$ and, consequently, 
\begin{equation*}
    e^{i\tfrac{2\pi s(x)}{N}}\cdot e^{i\tfrac{2\pi s(y)}{N}} = e^{i\tfrac{2\pi s(z)}{N}}.
\end{equation*}

Let us now consider the right hand-side of the \eqref{eq:group_bijective}. It is equivalent to the following: there exists $k \in \Z$ such that $s(x) + s(y) = s(z) + kN$. Since $s(x), s(y), s(z) \in \chi$, then either $k = 0$ or $k = 1$, which is equivalent to the left hand-side of the \eqref{eq:group_bijective}. $\Box$

Theorem \ref{th:operation} also shows that for random variable $X$ taking values in $\chi$ with pseudo-summation \eqref{eq:main_lut} there exists a classic characteristic function. Indeed, function 
\begin{equation*}
    f_X(t) = \sum_{k=0}^{N-1} \PP(X = k) e^{i\tfrac{2\pi s(k)}{N}t}, \quad t = 1,\ldots, N-1
\end{equation*}
uniquely defines the random variable $X$ and if $X$ and $Y$ are independent, then 
\begin{equation*}
    f_{X\oplus Y}(t) = f_X(t)\cdot f_Y(t), \quad t = 1,\ldots, N-1.
\end{equation*}

In terms of matrices $U$ and $V$ from \eqref{eq:general_ch_fun} we have the following:
\begin{equation*}
    U = V = \Bigr( e^{i\tfrac{2\pi s(k)}{N}t} \Bigr), k,t = 0, \ldots, N-1. 
\end{equation*}

Without loss of generality, we can consider only identical permutation. Indeed, characteristic function of $Y = s^{-1}(X)$ equals the following:
\begin{equation*}
    f_Y(t) = \sum_{k=0}^{N-1} \PP\bigr(Y = k\bigr) e^{i\tfrac{2\pi s(k)}{N}t} = \sum_{k=0}^{N-1} \PP\bigr(X = s(k) \bigr) e^{i\tfrac{2\pi s(k)}{N}t} = \sum_{l=0}^{N-1} \PP\bigr(X = l \bigr) e^{i\tfrac{2\pi l}{N}t}, \quad t = 1,\ldots, N-1. 
\end{equation*}

Further, all the results will be formulated only for the identical permutation. The rest of the results are obtained by replacing the variable $Y = s(X)$. 

Similarly to the general case, we first consider stable distributions. To highlight this particular case, we will write $\Z_N$-stable random variables instead of $A$-stable ones.
\begin{Theorem}\label{th:Z-N-stable}
    Random variable $\xi$ is $\Z_N$-stable if and only if it is either degenerate at point $0$ or there exists a divisor $M$ of a number $N$ such that $\xi \disteq M\cdot U$, where $U$ is a uniformly distributed random variable on the set $[0:\tfrac{N}{M}]$.  
\end{Theorem}
\textbf{Proof.} \quad Random variable $\xi$ is $\Z_N$-stable if and only if for every $t = 0,1,\ldots, N-1$ we have
\begin{equation}\label{eq:stable}
    f^2_{\xi}(t) = f_{\xi}(t).
\end{equation}

First of all, let us notice that a degenerate law at point $0$ satisfy this condition. If $M$ is a divisor of a number $N$ and $U$ is a uniformly distributed random variable on the set $[0:R]$, where $R = N/M$. Then for all $t = z\cdot R + \tau$, $z = 0,\ldots, M-1$, $\tau = 1, \ldots, R-1$ the characteristic function of random variable $\xi = M\cdot U$ equals
\begin{equation*}
    f_{\xi}(t) = \sum_{k=0}^{R-1}\frac1{R} e^{i\tfrac{2\pi Mk}{N} (zR + \tau)} = \frac1{R} \sum_{k=0}^{R-1}e^{i\tfrac{2\pi k}{R} (zR + \tau)} = \frac1{R} \sum_{k=0}^{R-1}e^{i\tfrac{2\pi k}{R}\tau}.
\end{equation*}

If $\tau \ne 0$, then $f_{\xi}(t) = 0$ and $f_{\xi}(t) = 1$, otherwise (see \cite[Theorem 6.10]{Apostol}).

Let us show the opposite. Form \eqref{eq:stable} it follows that for every $t = 1,\ldots, N~-~1$ the function $f_{\xi}(t)$ is either $0$ or $1$. If $f_{\xi}(1) = 1$, then $\xi$ can only be degenerate law at point $0$. 

Now let us assume that for some $t = 2,\ldots, N-1$ function $f_{\xi}(t) = 1$. Let $t_* \in \{2,\ldots, N-1\}$ be such a number that $f_{\xi}(t_*) = 1$ and for all $0 < t < t_*$ the function $f_{\xi}(t) = 0$. Then 
\begin{equation}\label{eq:inequality}
    1 = \sum_{k=0}^{N-1} p_k e^{i\tfrac{2\pi k}{N} t_*} = \textrm{Re}\, \Bigr( \sum_{k=0}^{N-1} p_k e^{i\tfrac{2\pi k}{N} t_*} \Bigr) = \sum_{k=0}^{N-1} p_k \cos\bigr( \tfrac{2\pi k}{N} t_*\bigr) \leqslant 1.
\end{equation}

For \eqref{eq:inequality} to be equality it is necessary and sufficient that $\PP(X\cdot t_* \,\%\, N = 0) = 1$, i.e.
\begin{equation*}
    \sum_{k\colon\, k\cdot t_*\,\%\,N = 0} p_k = 1.
\end{equation*}

It follows that, first of all, $t_*$ is a divisor of a number $N$ and, secondly, there exists a random variable $Y$ that takes values in the set $\{0, 1, \ldots, N/t_*-1\}$ such that $\xi \disteq t_*\cdot Y$. Since for all $0<t <t_*$ characteristic function $f_{\xi}(t) = 0$, then for all $\tau = 1,\ldots, N/t_* - 1$ characteristic function $f_{Y}(\tau) = 0$. It is \cite[Theorem 6.10]{Apostol} that states that $Y$ has uniform distribution on $[0: R]$ and concludes the proof. $\Box$

Let us note that Theorem \ref{th:Z-N-stable} possibly can be proved earlier. However, the authors failed to find the original proof. 

\begin{Corollary}
    If $N$ is a prime number, then random variable $\xi$ is $\Z_N$-stable if and only if it is either degenerate at point $0$ or uniformly distributed on $\{0,1,\ldots, N-1\}$. 
\end{Corollary}

The following theorem provides the necessary and sufficient conditions for random variable $X$ to belong to the domain of attraction of each stable law.
\begin{Theorem}\label{th:Z_stable}
    Let $M$ be a divisor of a number $N$, $R = \tfrac{N}{M}$, and random variable $U$ has uniform distribution on the set $[0:R]$. Random variable $X$ belongs to the domain of attraction of $\xi = M\cdot U$ if and only if there exists random variable $Y$ such that $X = M\cdot Y$ and for all $a \in [0:R]$ and $t_* = 1,\ldots, R-1$ we have
\begin{equation*}
    \PP\Bigr( \bigr(Y - a\bigr)t_* \,\%\, R = 0\Bigr) < 1.
\end{equation*}
    
\end{Theorem}
\textbf{Proof.} \quad It has already been shown that characteristic function of $\xi$ equals 
\begin{equation*}
    f_{\xi}(t) = \begin{cases}
        1, & t\,\%\, M = 0; \\ 
        0, & t\,\%\, M \ne 0.
    \end{cases}
\end{equation*}

It follows that $X$ belongs to the domain of attraction of $\xi$ if and only if for every $t = 1,\ldots, N-1$ we have
\begin{equation*}
    \lim\limits_{m\to\infty} \Bigr(f_{X}(t)\Bigr)^m = \begin{cases}
        1, & t\,\%\, M = 0; \\ 
        0, & t\,\%\, M \ne 0.
    \end{cases}
\end{equation*}

In particular, $\lim\limits_{m\to\infty} \Bigr(f_{X}(M)\Bigr)^m = 1$ and, as a consequence, $f_X(M) = 1$. It has already been shown at the proof of the Theorem \ref{th:Z_stable} that there exists a random variable $Y$ taking values in $[0, \tfrac{N}{M}]$ such that $X \disteq M\cdot Y$. 

Let $\{X_k\}$, $k \in \N$ be the sequence of independent copies of $X$, then 
\begin{equation*}
    \bigoplus_{k=1}^m X_k \disto M\cdot U, \quad m \to \infty 
\end{equation*}
if and only if 
\begin{equation*}
    \Bigr( \sum_{j=1}^m Y_j \Bigr) \, \%\, \frac{N}{M} \disto U, \quad m \to \infty.
\end{equation*}

Now, without loss of generality we can assume that $M = 1$. Then a random variable $X$ belongs to the domain of attraction of random variable $U$ with uniform distribution on $[0:N]$ if and only if for all $t = 1,\ldots, N-1$ we have
\begin{equation}\label{eq:limit}
    \lim\limits_{m\to\infty} \Bigr(f_{X}(t)\Bigr)^m = 0.
\end{equation}

One can see that \eqref{eq:limit} holds if and only if for all $t = 1,\ldots, N-1$ we have $|f_X(t)| < 1$ . Let assume that there exists $t_* = 1,\ldots, N-1$ and $a \in \{0,\ldots, N-1\}$ such that 
\begin{equation*}
    f_X(t_*) = \sum_{k=0}^{N-1} \PP\bigr(X = k\bigr) e^{i\tfrac{2\pi k}{N}t_*} = e^{i\tfrac{2\pi a}{N}t_*}.
\end{equation*}

It follows that $X$ does not belong to the domain of attraction of $U$ if and only if there exists $a \in \{0,\ldots, N-1\}$ and $t_* = 1,\ldots, N-1$ such that 
\begin{equation*}
    \PP\Bigr( \bigr(X - a\bigr)t_* \,\%\, N = 0\Bigr) = 1. \quad \Box
\end{equation*}

As noted earlier, $\Bigr(\chi, \oplus\Bigr)$ is a group, moreover, it is locally compact and Abelian. For such groups, a complete classification of infinitely divisible distributions in terms of their characteristic functions is known (see \cite{Rao}). Recall that a random variable $\xi$ is called \textit{infinitely divisible} if for every $n \in \N$ there exists a random variable $X^{(n)}$  such that
\begin{equation*}
    \xi \disteq X_1^{(n)} \oplus \ldots \oplus X_n^{(n)} \oplus a_n,
\end{equation*}
where $X_1^{(n)}, \ldots, X_n^{(n)}$ are independent copies of $X^{(n)}$, $a_n \in \chi$. 

\begin{Theorem}\label{th:id}
    For any infinitely divisible random variable $\xi$ there exists $a \in \chi$, a divisor $M$ of $N$, and a random variable $X$ such that 
    \begin{equation*}
        \xi \disteq a \oplus \bigr(M\cdot U) \oplus \bigoplus_{k=1}^{\mathcal{N}} X_j,
    \end{equation*}
    where $U$ has a uniform distribution  on the set $[0, \tfrac{N}{M}]$, $\{X_k\}$, $k \in \N$ are independent copies of $X$, and $\mathcal{N}$ has the Poisson distribution with intensity $\lambda \geqslant 0$.   
\end{Theorem}
\textbf{Proof.} \quad From \cite[Theorem 7.1]{Rao} the random variable $\xi$ is infinity divisible if and only if its characteristic function has the following representation
\begin{equation*}
    f_{\xi}(t) = \lambda_0(t) \cdot \exp\Bigr\{ i\frac{2\pi a}{N}t - \varphi(t) + \sum_{k=0}^{N-1} c_k\bigr( e^{i\tfrac{2\pi k}{N}t} - 1 - i\tfrac{2\pi k}{N}t \bigr) \Bigr\},
\end{equation*}
where $\lambda_0(t)$ is a characteristic function of $M\cdot U$ with some $M$, $c_k \geqslant 0$, $\sum_{k=0}^{N-1}c_k < \infty$, for all $t_1,t_2 \in \chi$ the function $\varphi$ satisfy the following condition:
\begin{equation}\label{eq:norm}
    \varphi\Bigr( \bigr(t_1 + t_2\bigr) \,\%N \Bigr) + \varphi\Bigr( \bigr(t_1 - t_2\bigr) \,\%N \Bigr) = 2\varphi(t_1) + 2\varphi(t_2). 
\end{equation}

To prove the Theorem \ref{th:id} it is sufficient to show that $\varphi(t) = 0$ for all $t \in \chi$. Indeed, $\lambda_0(t)$ corresponds to the random variable $M\cdot U$, $i\frac{2\pi }{N}t\bigr(a - \sum_{k=0}^{N-1} c_k k\bigr)$ corresponds to the shift, and finally, if $\lambda = \sum_{k=0}^{N-1}c_k$ and $\PP( X_1 = k) = c_k/\lambda$, $k=0,\ldots, N-1$, then 
\begin{equation*}
    f_Y(t) = \E e^{i\tfrac{2\pi Y}{N}t} = \sum_{l=0}^{\infty}\frac{\lambda^l}{l!}e^{-\lambda} \Bigr( \E e^{i\tfrac{2\pi X_1}{N}t} \Bigr)^l = \exp\Bigr\{ \lambda \bigr(\E e^{i\tfrac{2\pi X_1}{N}t} -1 \bigr) \Bigr\}, \quad Y = \bigoplus_{k=1}^{\mathcal{N}} X_j. 
\end{equation*}

Since $\lambda = \sum_{k=0}^{N-1}c_k$ and $\PP( X_1 = k) = c_k/\lambda$, then 
\begin{equation*}
    f_Y(t) = \exp\Bigr\{ \lambda \sum_{k=0}^{N-1}\tfrac{c_k}{\lambda}\bigr(e^{i\tfrac{2\pi k}{N}t} -1 \bigr) \Bigr\} = \exp\Bigr\{\sum_{k=0}^{N-1} c_k\bigr(e^{i\tfrac{2\pi k}{N}t} -1 \bigr) \Bigr\}.
\end{equation*}

Now let us show that $\varphi(t) = 0$ for all $t =0,\ldots, N-1$. In \eqref{eq:norm} let suppose that $t_2 = 0$, then $2\varphi(t_1) = 2\varphi(t_1) + 2\varphi(0)$. It follows that $\varphi(0) = 0$. Let us prove that for every $t = 0,\ldots, N-1$ there exists an integer $k = k(t)$ such that $\varphi(t) = k\varphi(1)$ and, moreover, $k(t+1) > k(t)$. For $t = 0$ and $t = 1$ it is obvious that $k(0) = 0$ and $k(1) = 1$. 

In \eqref{eq:norm} let suppose that $t_1 = t$, $t_2 = 1$, $1 \leqslant t < N-1$ then 
\begin{equation*}
    \varphi(t+1) + \varphi(t-1) = 2\varphi(t) + 2\varphi(1).
\end{equation*}

From the induction step, it follows that 
\begin{equation*}
    \varphi(t+1) =  2\varphi(t) + 2\varphi(1) - \varphi(t-1) = \bigr( 2k(t) - k(t-1) + 2 \bigr)\varphi(1).
\end{equation*}

One can see that $k(t+1) = 2k(t) - k(t-1) + 2$ and from the induction step we have $k(t) - k(t-1) > 0$, which leads to the following inequality $k(t+1) > k(t) + 2 > k(t)$. 

Hence, it has been proved that for any $t < N-1$ there exists an integer $k(t)$ such that $\varphi(t) = k(t)\varphi(1)$ and $k(t-1) < k(t)$. 

In \eqref{eq:norm} let $t_1 = t_2$. Then 
\begin{equation}\label{eq:similar}
    \varphi\bigr( (2t_1)\,\%N \bigr) + \varphi(0) = 4\varphi(t_1). 
\end{equation}

Similarly, let $t_2 = N - t_1$. Then 
\begin{equation}\label{eq:diff}
     \varphi(0) + \varphi\bigr( (2t_1)\,\%N \bigr) = 2\varphi(t_1) + 2\varphi(N - t_1). 
\end{equation}

From \eqref{eq:similar} and \eqref{eq:diff} one can see that for any $t =1,\ldots, N-1$ we have $\varphi(t) = \varphi(N-t)$. If $\varphi(1) \ne 0$, then $k(t) = k(N-t)$. However, for $t < N/2$ it has been proved that $k(t) < k(N-t)$, which follows that $\varphi(1) = 0$ and, consequently, $\varphi(t) = 0$ for all $t =0,\ldots, N-1$. $\Box$

\textbf{B)} Another special case of a Look-Up Table is the maximum value, that is
\begin{equation*}
    x \oplus y = A(x,y) = \max\{x,y\}, \quad x,y \in \chi.
\end{equation*}

It is obvious that $\oplus$ is an associative and commutative operation. Also, for every $x \in \chi$ we have $x \oplus x = x$. In Section 2 it has been proved that degenerate law at point $x$ is stable for all $x \in \chi$. In this case, similar to the real-valued case, we will write that degenerate laws are max-stable. 

Let us show that there is no other max-stable laws. If $\xi$ is max-stable, then for every $x \in \chi$ we have
\begin{equation*}
    \PP\bigr( \xi_1 \oplus \xi_2 \leqslant x\bigr) = \PP\bigr(\xi \leqslant x\bigr),
\end{equation*}
where $\xi_1$, $\xi_2$ are independent copies of $\xi$. 

Since $x \oplus y = \max\{x,y\}$, then for every $x \in \chi$ we have 
\begin{equation*}
    \PP\bigr( \xi \leqslant x\bigr)^2 = \PP\bigr( \xi \leqslant x\bigr).
\end{equation*}

It means that for every $x \in \chi$ the probability $\PP\bigr( \xi \leqslant x\bigr)$ is either 0 or 1. Since $F$ is a non-decreasing function, then $\xi$ is degenerate random variable. These arguments entail the following theorem.
\begin{Theorem}
    Random variable $\xi$ is max-stable if and only if it is degenerate at some point $x \in \chi$. 
\end{Theorem}

Theorem \ref{th:dom_attr} yields the necessary and sufficient conditions for random variable $X$ to belong to the domain of attraction of a degenerative law at point $x$. Namely, the following theorem is true.
\begin{Theorem}
    Random variable $X$ belongs to the domain of attraction of a degenerate law at point $x$ if and only if $\PP(X > x) = 0$ and $\PP(X = x) > 0$. 
\end{Theorem}
\textbf{Proof.} \quad From Theorem \ref{th:dom_attr} it follows that if $X$ belongs to the domain of attraction of a degenerate law at point $x$, then $\PP(X > x) = 0$. If $\PP(X = x) = 0$, then $\PP(\max\{X_1,\ldots, X_m\} = x) = 0$ for all $m \in \N$, which entails a contradiction. 

Let now assume that $\PP(X > x) = 0$ and $\PP(X = x) > 0$. Then for $y < x$ we have
\begin{equation*}
    \PP(\max\{X_1,\ldots, X_m\} < y) = \Bigr( \PP(X < y) \Bigr)^m < \Bigr(1 - \PP(X = x)\Bigr)^m \to 0, \quad m \to \infty. 
\end{equation*}

Since $\PP(X > x) = 0$, then $\PP(\max\{X_1,\ldots, X_m\} > x) = 0$. It follows that $X_1 \oplus \ldots \oplus X_m$ converges to degenerative law at point $x$. $\Box$

Similar to the previous case, let us study infinitely divisible distributions. Let assume that random variable $\xi$ with the distribution function $F$ is infinitely divisible. It means that for every $n \in \N$ there exists a random variable $X^{(n)}$ with the distribution function $F_n$ such that
\begin{equation*}
    \xi \disteq X^{(n)}_1 \oplus \ldots \oplus X^{(n)}_n = \max\limits_{k=1,\ldots,n}\bigr\{ X^{(n)}_k\bigr\},
\end{equation*}
where $X^{(n)}_1, \ldots, X^{(n)}_n$ are independent copies of $X^{(n)}$. 

Thus, in terms of distribution functions, distribution function is infinitely divisible if and only if for every $n \in \N$ there exists a distribution function $F_n$ such that 
\begin{equation*}
    F^n_n(x) = F(x), \quad x \in \R. 
\end{equation*}

Let us fix an arbitrary distribution function $F$. It is piece-wise constant function with jumps in points $x \in \chi$, then for all $n \in \N$ the function $F_n(x) = \bigr(F(x) \bigr)^{1/n}$ is also a non-decreasing piece-wise constant function with jumps in points $x \in \chi$ and 
\begin{eqnarray*}
    \lim_{x\to-\infty}F_n(x) = 0, \,\, \lim_{x\to\infty}F_n(x) = 1. 
\end{eqnarray*}

It means that $F_n$ is a distribution function of a discrete law on $\chi$. Hence, if $x \oplus y = \max\{x,y\}$, then any random variable is infinitely divisible. 

Let us note that in this case Look-Up Table also satisfy \eqref{eq:general_ch_fun} with matrices
\begin{equation*}
    U = V = \Bigr( \mathbb{I}_{i \geqslant j} \Bigr), \, i,j = 0,\ldots, N-1,
\end{equation*}
where $\mathbb{I}_{i \geqslant j} = 1$ if $i \geqslant j$ and 0 otherwise. 

\section*{Acknowledgments}
The work of I. A. Alexeev was supported in part by the Moebius Contest Foundation for Young Scientists.

\end{document}